\documentclass[10pt]{amsart}

\usepackage{amssymb} 
\usepackage[latin1]{inputenc}
\usepackage[matrix,arrow,curve]{xy}
\usepackage[mathscr]{eucal}
\usepackage{epic,eepic}
\usepackage{bbm}
\usepackage{multirow}
\usepackage{graphicx}
% ...
\makeindex

\newtheorem{theorem}{Theorem}

\newtheorem{lemma}{Lemma}

\theoremstyle{definition}

\newtheorem{example}{Example}

\newenvironment{pf}
{\medskip\noindent {\it Proof. \ }}
{\hfill\nobreak $\Box$ \par\bigbreak}

\newtheoremstyle{numero}
{1ex}
{1ex}
{\it}
{1cm}
{\scshape}
{}
{4pt}
{}
\theoremstyle{numero}

\newcommand{\F}{ \mathbb F}

\newcommand{\Q}{{ \mathbb Q } }
\newcommand{\Z}{{ \mathbb Z  }}
\newcommand{\N}{{ \mathbb N  }}

\newcommand{\anneau}{{ \mathcal O}}

\newcommand{\End}{{\text{End}}}

\newcommand{\SL}{{\text {SL}}}

\newcommand{\Pc}{{\mathcal{P}}}

\newcommand{\G}{{\mathfrak g}}

\newcommand{\Ac}{{\mathcal{A}}}

\newcommand{\Frob}{{\rm{Frob}}}

\newcommand{\ord}{{\text{ord}}}

\renewcommand{\G}{G}

\begin{document}

\baselineskip15pt

\title{Non-zero coefficients of half-integral weight modular forms mod $\ell$} 
\author{Jo\"el Bella\"{\i}che, Ben Green, Kannan Soundararajan}
\thanks{Jo\"el Bella\"{\i}che was supported by NSF grant DMS 1405993.   Ben Green was supported by a Simons Investigator grant from the Simons Foundation.  Kannan Soundararajan was partially supported by NSF grant DMS 1500237, and a Simons Investigator grant from the Simons Foundation.   Part of the work was carried out when the second and third authors were in residence at MSRI, Berkeley during 
the Spring semester of 2017, supported in part by NSF grant DMS 1440140} 
\address{Department of Mathematics, Brandeis University, Waltham, MA 02453} 
\email{jbellaic@brandeis.edu}
\address{Mathematical Institute, University of Oxford, Oxford OX2 6GG} 
\email{ben.green@maths.ox.ac.uk}
\address{Department of Mathematics, Stanford University, Stanford, CA 94305} 
\email{ksound@stanford.edu} 
\maketitle 

\begin{abstract}
We obtain new lower bounds for the number of Fourier coefficients of a weakly holomorphic modular form of half-integral weight not divisible by some prime $\ell$. Among the applications of this we show that there are $\gg \sqrt{X}/\log \log X$ integers $n \leq X$ for which the partition function $p(n)$ is not divisible by $\ell$, and that there are $\gg \sqrt{X}/\log \log X$ values of $n \leq X$ for which $c(n)$, the $n$th Fourier coefficient of the $j$-invariant, is odd.
\end{abstract}

\section{Introduction}

\noindent Let $K$ be a number field and  $\anneau$ its ring of integers.  Let $\ell$ be a rational prime and let $\lambda$ be a maximal ideal of $\anneau$ above $\ell$. We denote by $\F$ the residue field $\anneau/\lambda$, a finite extension of $\F_\ell$.
 The reader will lose little by supposing that $K = \Q,$ $\anneau = \Z$, $\lambda = (\ell)$ and $\F=\F_\ell$; our main applications use only this case.

\begin{theorem}\label{thmform}  
Let $f=\sum_{n=n_0}^\infty a_n q^n$ be a weakly holomorphic modular form\footnote{
{\it Weakly holomorphic} allows for polar singularities at the cusps; for this and other basic definitions, we refer the reader to 
\cite[Chapter 1]{ono}.} 
of weight $k \in \frac{1}{2}\Z \setminus \Z$ and level $\Gamma_1(N)$.   Suppose that the coefficients $a_n$ lie in  the ring 
$\anneau$.   If $\ell \ge 3$, we assume that $f \not \equiv 0 \pmod{\lambda}$, and for $\ell =2$ we 
assume that $f \pmod{\lambda}$ is not a constant. Then
$$\# \{ n \leq X, \ a_n \not \equiv 0 \pmod{\lambda}\}  \gg \frac{\sqrt{X}}{\log \log X}.$$ 
\end{theorem}
\par \bigskip

Here are some sample applications of Theorem~\ref{thmform}.

\begin{example} Take $f=\eta_1(z)^{-1}$ with $\eta_1(z)=\eta(24 z)$ (Dedekind's eta function), so that $f$ is a weakly holomorphic modular form of weight $-1/2$ and level $\Gamma_0(576)$.   The Fourier expansion of $f$ is 
$$
f(q) = q^{-1} \prod_{n=1}^\infty \left({1-q^{24n}}\right)^{-1} = \sum_{n=0}^\infty p(n) q^{24n-1}
$$
where $p(n)$ is the {\it partition function} (cf. \cite[Theorem 1.60, Corollary 1.62 and Theorem 5.3]{ono} for a proof of these well-known facts). Applying the theorem to $f$, we conclude that 
\begin{equation} 
\label{1.1} 
\# \{ n \leq X, \ p(n) \not \equiv 0 \pmod{\ell}\}  
\gg \frac{\sqrt{X}}{\log \log X}.
\end{equation} 
This improves, by a factor of about $(\log X)^{\frac 34}$,  earlier results of Ahlgren \cite{ahlgren}, Chen \cite{chen} and Dai \& Fang \cite{DF} for odd $\ell$.  
In the case $\ell=2$, \eqref{1.1}  improves upon previous results (established by somewhat different methods than for $\ell$ odd) 
by a factor of about $(\log X)^{\frac 78}$; see    \cite{BN, nicolasRS, nicolasvi, nicolasParite}.
\end{example}

\begin{example} Theorem \ref{thmform} applies in particular when $f$ is a {\it holomorphic} cusp form of half-integral weight $k$ (which must then be positive).
In this case, it improves on the main theorem of ~\cite{BrOno} (itself an improvement of ~\cite{SkOno}) which proves the 
slightly weaker estimate $\# \{ n \leq X, a_n \not \equiv 0 \pmod{\ell}\}  \gg {\sqrt{X}}/{\log X}$ under the supplementary assumption that the coefficients $a_n \pmod{\lambda}$ are not supported in a finite union of sequences of the form $(c n^2)_{n \in \N}$.  We remark that in \cite{BrOno, SkOno} the Shimura  
correspondence between holomorphic half-integral weight modular forms and integral weight modular forms plays a crucial role, whereas 
our proof of Theorem \ref{thmform} does not involve the Shimura correspondence.  
\end{example}

\begin{example} When $\ell=2$, our theorem applies as well to weakly holomorphic modular forms of {\it integral} weight, since those forms
are congruent modulo $\lambda$ to forms of {\it half-integral} weight (see Lemma~\ref{weighthalf} below). In particular, for the modular invariant $j(q)=\sum_{n=-1}^\infty c(n) q^n$, which is of weight $0$ and level $\SL_2(\Z)$, we obtain that 
\begin{equation} 
\label{1.2} 
\# \{ n\le X, c(n) \text{  is odd} \} \gg \frac{\sqrt X}{\log \log X}.
\end{equation} 
This improves upon recent results in \cite{alfes, zanello} obtained by different methods.  

For completeness, we remark that \cite{BS} establishes, for $\ell \ge 3$, an asymptotic 
for  the number of non-zero coefficients (mod $\ell$) of holomorphic modular forms, and \cite{BN} establishes such an asymptotic for $\ell=2$ and holomorphic forms of level $1$.  
The situation for {\it weakly} holomorphic forms of integral weight remains mysterious, and  (for example) we do not have lower bounds for the number of $c(n) \neq 0 \pmod \ell$ 
for primes $3 \le \ell \le 11$.   
\end{example}

\par \bigskip

The proof of Theorem~\ref{thmform} uses the standard idea of multiplying $f$ by a suitable lacunary holomorphic cusp form $g$ of half-integral weight, to obtain an holomorphic cusp form $h=fg$ of integral weight.   The panoply of results stemming from the existence of  Galois representations associated to integral weight holomorphic eigenforms 
may then be used to study the coefficients of $h$. Finding a suitable form $g$ is easy when $\ell>2$,  and somewhat less so in the case $\ell=2$.  
Our improvement over previous work comes from analyzing more carefully the implications for non-vanishing coefficients of the equality $h=fg$, 
which leads to a problem in analytic number theory/additive combinatorics. 

\begin{theorem}\label{theoremAS}  Let $u \ge 1$ be a fixed natural number, and let $X$ be large.  For any subset $\Ac \subset \{1.\dots,X\}$ the 
number of primes $p$ such that $pu \le X$ and $pu$ may be written as $a + m^2$ for some $a \in \Ac$ and some integer $m$ is 
$$ 
\ll \frac{\sqrt{X}}{\log X} \Big( |\Ac| \log \log X+ |\Ac|^{\frac 12} X^{\frac 14} \Big). 
$$ 
 \end{theorem}

Our interest in Theorem~\ref{theoremAS} is in the situation where a positive proportion of the primes $p$ 
are known to be of the form $a+m^2$, when it follows that $|\Ac|$ must have $\gg \sqrt{X}/\log \log X$ elements.  
 This statement is optimal, as we shall show in section~\ref{optimality} by constructing an example of a set $\Ac$ with $|\Ac| \asymp   {\sqrt{X}}/{\log \log X}$, 
 and with a positive proportion of primes below $X$ being of the form $a+m^2$.

Theorem~\ref{thmform}, on the other hand, is almost certainly not optimal.  
For any weakly holomorphic form $f(q)=\sum a_n q^n$ of half-integral weight, one might expect 
 $$
 \# \{ n \leq X, \ a_n \not \equiv 0 \pmod{\lambda}\}  \gg \sqrt{X},
 $$ 
 and this bound is attained for $\eta_1(q)$ (see (\ref{eta1}) below). 
Theorem \ref{thmform} comes close to this estimate.  For most forms $f$ of half-integral weight however 
(specifically for $f(q)=\eta_1^{-1}(q)=\sum_n p(n) q^{24 n-1}$, and perhaps for all forms that are not congruent mod 
$\lambda$ to a one-variable theta series), it is expected  that $f \pmod {\lambda}$ is not {\it lacunary}, which 
is to say that $
\# \{ n \leq X, \ a_n \not \equiv 0 \pmod{\lambda}\}  \gg X$.  

\section{Deduction of Theorem~\ref{thmform} from Theorem~\ref{theoremAS}}

\subsection{A preliminary lemma} 

Let $M_k(\Gamma_1(N),\anneau)$ be the $\anneau$-module of holomorphic modular forms of integral weight $k\geq 0$, 
level $\Gamma_1(N)$, and coefficients in $\anneau$. Let $\anneau_\lambda$ be the completion of $\anneau$ at the place defined by the ideal $\lambda$, and set $M_k(\Gamma_1(N),\anneau_\lambda) = M_k(\Gamma_1(N),\anneau) \otimes_\anneau \anneau_\lambda$.
Let $A$ be the closure of the $\anneau_\lambda$-subalgebra of $\End_{\anneau_\lambda}(M_k(\Gamma_1(N),\anneau_\lambda))$ generated by the Hecke operators $T_n$ for $n$ running among integers relatively primes to $N \ell$. Denote by $G_{\Q,N\ell}$ the Galois group of the maximal extension of $\Q$ unramified outside $N\ell$, and by $\Frob_p$, for $p$ a prime not dividing $N\ell$, the Frobenius element of $p$ in $G_{\Q,N\ell}$, well-defined up to conjugation.

\begin{lemma} \label{ext}There exists a unique continuous map $t: G_{\Q,N\ell} \rightarrow A$ which is central and satisfies  $t(\Frob_p) = T_p$ for every prime $p$ not dividing $N\ell$. This maps also satisfies $t(1)=2$. 
\end{lemma}
\begin{pf} This follows from a well-known argument of Wiles based on the existence of Galois representations attached to eigenforms due to Deligne; see, for example, \cite[Thm 1.8.5]{B} for a detailed proof.
\end{pf}

\subsection{The case $\ell>2$} 

We begin with a lemma.

\begin{lemma}  \label{posdens}  Assume that $\ell$ is odd.   Let $k \in \N$ and $h(q) = \sum_{n=0}^\infty a_n q^n\in M_k(\Gamma_1(N),\anneau)$. 
Let $u \geq 1$ be an integer such that $a_u \not \equiv 0 \pmod{\lambda}$.  There is a positive density set of primes $\Pc$ such that $a_{up} \not \equiv 0 \pmod{\lambda}$
for every $p\in {\mathcal P}$. 
\end{lemma}
\begin{pf}
With $t$ as in Lemma~\ref{ext}, the map from $G_{\Q,N\ell}$ to $\F$ sending $g$ to $a_u( t(g) h) \pmod{\lambda}$  is a continuous map. Thus there exists an open neighborhood $U$ of $1$ in  $\G_{\Q,N\ell}$ such that  $g \mapsto a_u( t(g) h) \pmod{\lambda}$ is constant on $U$. Let $\Pc$ be the set of primes not dividing $N\ell u$ such that $\Frob_p \in U$. By Chebotarev, $\Pc$ has positive density.  Further for $p \in \Pc$ we have 
$$
a_{up}(h) = a_u(T_p h) = a_u(t(\Frob_p) h) =   a_u(t(1) h) = 2 a_u(h)     \not\equiv 0\pmod{\lambda},
$$ 
since $\ell \neq 2$ and $a_u(h) \not \equiv 0 \pmod{\lambda}$.
\end{pf}

We can now deduce Theorem~\ref{thmform} from Theorem~\ref{theoremAS}.  
Let $f = \sum_{n \geq n_0} a_n q^n$ be a weakly holomorphic modular form of half-integral weight, level $\Gamma_1(N)$, and  coefficients in $\anneau_\lambda$.
Let $\eta(z)$ be the usual Dedekind's eta function and set $\eta_1(z)=\eta(24 z)$ so that (see \cite{ono})
\begin{eqnarray}
 \label{eta1} \eta_1(q) = q \prod_n (1-q^{24n}) = \sum_{n=-\infty }^\infty (-1)^n q^{(6n+1)^2} 
 \end{eqnarray}
 is a holomorphic cuspidal modular form of weight $1/2$.  
   Let $m$ be an even integer such that $\ell^m$ is larger than the order of any pole of $f$. 
 Then $h = f \eta_1^{\ell^m}$ is a holomorphic cuspidal modular form of integral weight $k+\ell^m/2$. 
 Since $f\pmod \lambda$ and $\eta_1 \pmod{\lambda}$ are non-zero, the power series $h \pmod{\lambda} \in \F[[q]]$ is also non-zero, 
 and indeed $h \pmod{\lambda}$ is not a constant (because $h$ is cuspidal, and a cuspidal constant form must be $0$). 

 Let $\Ac=\{n,\ a_n =a_n(f) \not \equiv 0 \pmod{\lambda}\}$.   Note that from \eqref{eta1}  
 $$ 
 \eta_1(q)^{\ell^m} = \Big( \sum_{n=-\infty}^{\infty} (-1)^{n} q^{(6n+1)^2} \Big)^{\ell^m} 
 \equiv \sum_{n=-\infty}^{\infty} (-1)^n q^{\ell^m (6n+1)^2} \pmod \lambda,  
 $$ 
 so that the Fourier coefficients of $\eta_1^{\ell^m}$ are non-zero $\pmod m$ only on squares.   Thus if $n$ is 
 such that $a_n(h) \not \equiv 0 \pmod \lambda$, then $n$ must be of the form $a + m^2$ for some $a\in \Ac$ and 
 some integer $m$.   Now, by Lemma~\ref{posdens}, the set of $n$ such that $a_n(h) \not \equiv 0 \pmod{\lambda}$ contains a set of the form 
 $u \Pc$, for some fixed natural number $u$ and a set of primes $\Pc$   of positive density.   
  For large $X$, it follows from Theorem \ref{theoremAS} that the number of primes $p$ with $up \le X$ and $up$ of the form 
  $a+ m^2$ with $a\in \Ac$ is $\ll |\Ac \cap [1,X]| \sqrt{X} (\log \log X)/\log X + |\Ac\cap [1,X]|^{\frac 12} X^{\frac 34}/\log X$.   It follows that 
  $|\Ac \cap [1,X]| \gg {\sqrt{X}}/{\log \log X}$, proving Theorem \ref{thmform}.

\subsection{The case $\ell=2$.} This case needs a little more care, and we begin by recalling a well-known result that (for $\ell =2$) modular forms of integer weights are 
congruent to modular forms of half-integer weight.  

\begin{lemma} \label{weighthalf} Assume that $\ell=2$.    For every weakly holomorphic modular form $f$ of weight $k$ and level $\Gamma_1(N)$ with coefficients in $\anneau$ there exists a weakly holomorphic modular form $f'$ of weight $k+1/2$, some level $\Gamma_1(N')$, with coefficients in $\anneau$, such that $f \equiv f' \pmod{\lambda}$.
\end{lemma}
\begin{pf} Recall that (see, for example, \cite[Prop 1.4]{ono}) the theta series $\theta_0(q)=\sum_{n=-\infty}^\infty q^{n^2} = 1 + \sum_{n=1}^\infty 2 q^{n^2}$ is a holomorphic modular form of weight $1/2$, level $\Gamma_0(4)$, coefficients in $\Z$.   Now take $f'=f \theta_0$.
\end{pf}

\begin{lemma} \label{uglylemma} Let $n_0$ be a non-zero integer, and let $N$ be a positive integer. There are only finitely many natural numbers $m$ such that 
$2^m + n_0$  equals $uy^2$ for some square-free divisor $u$ of $2N$.   
\end{lemma}
\begin{pf}
Write $m= 3 m_0 + r$ with $r=0$, $1$, or $2$, and set $x= 2^{m_0}$.  The equation $2^m+n_0 = uy^2$ becomes $2^r x^3 +n_0 = uy^2$, 
which for a given $r$ and $u$ may be viewed as an elliptic curve (since $n_0 \neq 0$).   By Siegel's theorem there are only finitely many 
integer points $(x,y)$ on this elliptic curve.   Since there are only three possible values for $r$, and finitely many possibilities for $u$ (being a 
square-free divisor of $2N$), the lemma follows. 
 \end{pf}

\begin{lemma} \label{uP1} Let $k \in \N$ and $h(q) = \sum_{n=0}^\infty a_n q^n\in M_k(\Gamma_1(N),\anneau)$. 
 Assume that there exists an integer $n \geq 1$ and a prime $p_0$ not dividing $2N$ such that $\ord_{p_0}{n}$ is odd and
$a_{n} \not \equiv 0 \pmod{\lambda}$. 
Then there exists an integer $u \geq 1$ and a set of primes $\Pc$ of positive density such that $a_{up} \not\equiv 0 \pmod{\lambda}$ for every $p \in \Pc$.
\end{lemma}
\begin{pf} 
We recall that if $h =\sum_{n=0}^\infty a_n q^n$ is a modular form for $\Gamma_1(N)$, and if $p$ is a prime not dividing $N$,
the $m$-th coefficient of the form $T_p h$ is
\begin{eqnarray} \label{hecke} a_m(T_p h) = a_{mp}(h) + p^{k-1} a_{m/p}(\langle p \rangle h),\end{eqnarray}
where $\langle p \rangle$ is the diamond operator, and with the convention that $a_{m/p}(-) = 0$ when $p \nmid m$. 

We claim that {\it if $p$ is a prime not dividing $2N$ and if $T_p h \pmod{\lambda}$ is a constant, then $a_n(h)=0$ if $\ord_p(n)$ is odd}.
We prove that claim, for all $h$ such that $T_p h \pmod{\lambda}$ is constant, by induction over the odd number $\ord_{p}(n)$. If $\ord_{p}(n)=1$,
applying (\ref{hecke}) to the form $h$ and the integer $m=n/p$ and reducing mod $\lambda$ gives (using that $p \equiv -1 \equiv 1 \pmod{\lambda}$):
$$a_{n}(h) \equiv a_{n/p^2} (\langle p \rangle h) = 0 \pmod{\lambda}.$$
For a general $n$ with $\ord_p(n)$ odd, we get similarly
$$a_{n}(h) \equiv a_{n/p^2} (\langle p \rangle h)\pmod{\lambda}.$$
By the induction hypothesis applied to the form $\langle p \rangle h$ (which also satisfies $T_p (\langle p \rangle h) \pmod{\lambda}$ constant since the diamond operator $\langle p \rangle$ commutes with $T_p$ and stabilizes the subspace of constants),
we get $a_n(h) \equiv 0 \pmod{\lambda}$ which completes the induction step.

By the hypothesis of the Lemma, it follows that $T_{p_0} h \pmod{\lambda}$ is not a constant, that is to say 
there exists $u \geq 1$ such that $a_u (T_{p_0} h) \not \equiv 0 \pmod{\lambda}$, or equivalently, with $t$ as in Lemma~\ref{ext},
$a_u(t(\Frob_{p_0}) h) \not \equiv 0 \pmod{\lambda}$. By continuity of $t$, there exists an open set $U$ in $G_{\Q,Np}$ such that for $p$ a prime not dividing $Nu$, if $\Frob_p \in U$, then 
 $$
 a_{up}(h) = a_u(T_p h) = a_u(t(\Frob_p) h) \equiv a_u(t(\Frob_{p_0}) h) \not \equiv 0 \pmod{\lambda}.
 $$ The set $\Pc$ of such primes $p$ is a set of primes of positive density by Chebotarev.
\end{pf}
 
 We are now ready to prove Theorem~\ref{thmform} in the case $l=2$ using Theorem~\ref{theoremAS}. Let $f$ be a weakly holomorphic modular form of half-integral weight with coefficients in $\anneau_\lambda$. By Lemma~\ref{weighthalf} we may assume that $f$ has integral weight instead.
 
We consider the form $h:=f \eta_1^{2^m}$ for a suitable $m$ that will be specified below. We observe that if $m \geq 1$,
$\eta_1^{2^m}$ has integral weight, and so does $h$. Moreover, since $\eta_1$ is cuspidal, $h$ is also cuspidal holomorphic for $m$ large enough. 

Write $f \equiv \sum_{n = n_0}^\infty a_n q^n \pmod{\lambda}$ with $a_{n_0} \not \equiv 0$. 
The first term of $h \pmod{\lambda}$ is $a_{n_0} q^{2^m+n_0}$. When $n_0 \neq 0$, Lemma~\ref{uglylemma} ensures that we can choose $m$ even, large enough in the sense of the preceding paragraph, and such that there is a prime $p_0$ not dividing $2N$ such that $\ord_{p_0}(2^m+n_0)$ is odd. When $n_0=0$, let $a_{n_1} q^{n_1}$ with $n_1>0$, $a_{n_1} \not \equiv 0$ the term of smallest positive degree in $f \pmod{\lambda}$ (such an $n_1$ exists 
because we assume that $f \pmod{\lambda}$ is not a constant). If $m$ is such that $2^m > n_1$, then the form $h$ has a term $a_{n_1} q^{2^m+n_1}$. 
Again by Lemma~\ref{uglylemma} we can find $m$ large enough and such that there is a prime $p_0$ not dividing $2m$ such that $\ord_{p_0}(2^m+n_1)$ is odd.

So in both cases ($n_0 \neq 0$ and $n_0 =0$) we have shown the existence of an integer $m$ such that $h =f \eta_1^{2^m}$ is a cuspidal holomorphic modular form of integral weight and such that, by Lemma~\ref{uP1}, there is $u\ge 1$ and a set of primes $\Pc$ of positive density, with $a_{u p}(h) \not \equiv 0 \pmod{\lambda}$ for every $p \in \Pc$. The rest of the proof is now exactly as in the case $\ell>2$.

\section{Proof of Theorem \ref{theoremAS}}

\noindent Given a positive integer $a$, we let $\chi_{-4a}=(\frac{-4a}{\cdot})$ denote the Kronecker symbol, 
which is a Dirichlet character $\pmod{4a}$.   Note that $-4a$  is a discriminant, but it need not be a fundamental discriminant.  We 
denote the associated (negative) fundamental discriminant by ${\widetilde a}$, so that $-4a   = {\widetilde a} a_2^2$ for a suitable 
natural number $a_2$.   
 
\begin{lemma} \label{lemmalargesieve} Let $u$ be a fixed natural number, and let $X$ be large.  For every integer $1\le a \leq X$, 
$$
\# \{ p : \  up \le X, \  p = a+ m^2 \textrm{  for some } m\in {\Bbb Z} \} \ll  \frac{\sqrt{X}}{\log X} \prod_{p\le X^{\frac 14}} \Big( 1- \frac{\chi_{-4a}(p)}{p} \Big).
$$
\end{lemma}
\begin{pf}   Consider the equivalent problem of estimating the number of $m$ below $\sqrt{X}$ such that $a+ m^2$ is of the 
form $ur$ for a prime number $r$.   We may restrict attention to $r> X^{\frac 14}$, since the smaller primes $r$ contribute 
negligibly to the number of $m$.   For each prime $p \nmid 2u$ and $p\le X^{\frac 14}$ we see that $m^2$ cannot be $\equiv -a \pmod{p}$, 
which means that $1+ \chi_{-4a}(p)$ residue classes $\pmod p$ are forbidden for $m$.  Any standard upper bound sieve (for example, Brun's sieve 
or Selberg's sieve; or see Theorem 2.2 of \cite{HR}) then shows that the number of possible $m \le \sqrt{X}$ is 
$$ 
\ll \sqrt{X} \prod_{\substack { p\le X^{\frac 14} \\ p\nmid 2u}} \Big( 1- \frac{1+\chi_{-4a}(p)}{p} \Big) 
\ll \frac{\sqrt{X}}{\log X} \prod_{ p\le X^{\frac 14}} \Big( 1- \frac{\chi_{-4a}(p)}{p} \Big),
$$  
and the lemma follows.  
  \end{pf}
  
 Call a fundamental discriminant $d$ {\it good} if the corresponding Dirichlet $L$-function $L(s,\chi_d)$ 
 has no zeros in the region $\{ \sigma > 99/100, \ |t| \le |d| \}$, and call the discriminant $d$ {\it bad} otherwise. 


\begin{lemma} \label{agood}  Suppose $1\le a \le X$ is an integer, and that the  fundamental discriminant ${\widetilde a}$ corresponding to $-4a$ is good.   
Then 
$$
\prod_{p\le X^{\frac 14}} \Big( 1- \frac{ \chi_{-4a}(p)}{p} \Big) \ll {\log \log X}. 
$$ 
\end{lemma}
\begin{pf} 
By \cite[Lemma 2.1]{GS}, for a good fundamental discriminant ${\widetilde a}$ one has 
$$
L\Big(1+ \frac{1}{\log X}, \chi_{\widetilde a}\Big) \asymp \prod_{p < (\log |\widetilde a|)^{100}} \Big(1 - \frac{\chi_{\widetilde a}(p)}{p} \Big)^{-1}.
$$
Further 
\begin{align*}
\log L\Big(1+ \frac{1}{\log X}, \chi_{\widetilde a}\Big) &= \sum_{p} \frac{\chi_{\widetilde a}(p)}{p^{1+1/\log X} } + O(1)  \\
&= \sum_{p \le X^{\frac 14}} \frac{\chi_{\widetilde a}(p)}{p} + O\Big( \sum_{p\le X^{\frac 14}} \frac{1-p^{-1/\log X}}{p}  + \sum_{p > X^{\frac 14}} \frac{1}{p^{1+1/\log X}} +1\Big). 
\end{align*} 
Using $1-p^{-1/\log X} =O(\frac{\log p}{\log X})$ for $p\le X^{\frac 14}$, the first error term above is seen to be $O(1)$, and 
partial summation shows that the second term is also $O(1)$.   Therefore 
\begin{equation} 
\label{3.1}
\prod_{p\le X^{\frac 14}} \Big( 1- \frac{\chi_{\widetilde a}(p)}{p} \Big) \asymp L(1+1/\log X, \chi_{\widetilde a})^{-1} 
\asymp \prod_{p \le (\log |\widetilde a|)^{100}} \Big(1 - \frac{\chi_{\widetilde a}(p)}{p} \Big).
\end{equation}

Now write $-4a = {\widetilde a} a_2^2$ for some positive integer $a_2 \le \sqrt{X}$.  Then 
$$ 
\prod_{p\le X^{\frac 14}} \Big( 1- \frac{\chi_{-4a}(p)}{p} \Big) = 
\prod_{p\le X^{\frac 14} } \Big( 1- \frac{\chi_{\widetilde a}(p)}{p} \Big) \prod_{\substack{p\le X^{\frac 14} \\ p|a_2}} \Big( 1- \frac{\chi_{\widetilde a}(p)}{p} \Big)^{-1},  
$$
and using \eqref{3.1} this is 
\begin{equation} 
\label{3.2}
\asymp \prod_{\substack{ p\le (\log |{\widetilde a}|)^{100} \\ p\nmid a_2 }} \Big(1 - \frac{\chi_{\widetilde a}(p)}{p} \Big) 
\prod_{ \substack{ X^{\frac 14} \ge p \ge (\log |{\widetilde a}|)^{100} \\ p|a_2} }  \Big(1 - \frac{\chi_{\widetilde a}(p)}{p} \Big)^{-1}. 
\end{equation} 

The first product in \eqref{3.2} is clearly at most
$$ 
\prod_{p \le (\log |{\widetilde a}|)^{100}} \Big(1 + \frac 1p \Big) \ll \log \log |{\widetilde a}|. 
$$ 
As for the second product in \eqref{3.2}, this is 
$$ 
\le \prod_{ \substack{ X^{\frac 14} \ge p \ge (\log |{\widetilde a}|)^{100} \\ p|a_2} }  \Big(1 - \frac{1}{p} \Big)^{-1} 
\le \prod_{(\log {\widetilde a})^{100} \le p \le (\log {\widetilde a})^{100} + (\log X)^2} \Big(1-\frac{1}{p} \Big)^{-1}, 
$$ 
 since $a_2$ has at most $\log X$ prime factors, and  the product is largest if these prime factors are the first $\le \log X$ 
 primes all larger than $(\log {\widetilde a})^{100}$.  This quantity is easily seen to be $\ll \max(1, \frac{\log \log X}{\log \log |{\widetilde a}|})$, 
 proving the lemma.
 \end{pf}

Applying Lemmas \ref{lemmalargesieve} and \ref{agood} we see that the number of primes $p$ 
with $up \le X$ and $p$ of the form $a + m^2$ with $a \in {\Ac}$ coming from a {\it good} associated fundamental 
discriminant ${\widetilde a}$ is bounded by 
$$ 
\sum_{\substack{a \in {\Ac} \\ {\widetilde a} \text{  good} }} \frac{\sqrt{X}}{\log X} \prod_{p\le X^{\frac 14}} \Big( 1- \frac{\chi_{-4a}(p)}{p} \Big) 
\ll |{\Ac}| \frac{\sqrt{X}}{\log X} \log \log X. 
$$ 
It remains to bound the number of primes arising from {\it bad} fundamental discriminants ${\widetilde a}$.   Note that, with $-4a ={\widetilde a} a_2^2$,  
\begin{align*} 
\prod_{p\le X^{\frac 14}} \Big(1- \frac{\chi_{-4a}(p)}{p} \Big) 
&\ll \prod_{p\le X^{\frac 14}} \Big(1-\frac{\chi_{\widetilde a}(p)}{p} \Big) \prod_{p|a_2} \Big(1-\frac 1p\Big)^{-1} 
\\
&\ll \frac{a_2}{\phi(a_2)} L(1+1/\log X, \chi_{\widetilde a})^{-1} \ll |{\widetilde a}|^{\epsilon} \frac{a_2}{\phi(a_2)}, 
\end{align*}
where the final estimate follows by an obvious modification to Siegel's theorem which gives $L(1+1/\log X, \chi_{\widetilde a}) \gg |{\widetilde a}|^{-\epsilon}$.  
 Thus the number of primes arising from bad fundamental discriminants is 
 $$ 
  \ll \frac{\sqrt{X}}{\log X} \sum_{ \substack{{a\le X} \\ {{\widetilde a} \ \text{ bad}}}} \prod_{p\le X^{\frac 14}} \Big(1- \frac{\chi_{-4a}(p)}{p} \Big) 
  \ll \frac{\sqrt{X}}{\log X}\sum_{ \substack{{|{\widetilde a}| \le X} \\ {{\widetilde a} \ \text{ bad}}}} |\widetilde a|^{\epsilon} \sum_{\substack{-4a ={\widetilde a}a_2^2 \\ a\in \Ac}} \frac{a_2}{\phi(a_2)}.
$$
For a given ${\widetilde a}$ we may bound the sum over $a_2$ above using Cauchy-Schwarz; thus 
$$ 
 \sum_{\substack{-4a ={\widetilde a}a_2^2 \\ a\in \Ac}} \frac{a_2}{\phi(a_2)} 
 \le \Big( \sum_{\substack{-4a ={\widetilde a}a_2^2 \\ a\in \Ac}} 1 \Big)^{\frac 12} \Big(  \sum_{a_2 \le \sqrt{X/|{\widetilde a}|}} \Big(\frac{a_2}{\phi(a_2)} \Big)^2 \Big)^{\frac 12} 
 \ll \sqrt{|\Ac|}  \frac{X^{\frac 14}}{|{\widetilde a}|^{\frac 14}}. 
 $$ 
 We conclude that the number of primes $p$ arising from bad fundamental discriminants is 
 \begin{equation} 
 \label{3.3}  
 \ll \frac{\sqrt{X}}{\log X} |\Ac |^{\frac 12} X^{\frac 14} \sum_{\substack{|{\widetilde a}| \le X \\ {\widetilde a} \text{ bad } } } \frac{|{\widetilde a}|^{\epsilon}}{|\widetilde a|^{\frac 14}} 
 \ll \frac{|\Ac |^{\frac 12} X^{\frac 34}}{\log X}  \sum_{\substack{|{\widetilde a}| \le X \\ {\widetilde a} \text{ bad } } } |{\widetilde a}|^{-\frac 16}, 
 \end{equation} 
 upon choosing $\epsilon =1/12$.   At this stage we note that bad fundamental discriminants are rare by a standard zero density result (see for example \cite[Theorem 20]{Bom}): 
 thus there are at most  $\ll Y^{1/10}$ bad fundamental discriminants $d$ with $Y \le |d| \le 2Y$.   Therefore the sum over bad ${\widetilde a}$ in \eqref{3.3} converges, and 
 we conclude that the quantity in \eqref{3.3} is $\ll |\Ac |^{\frac 12} X^{\frac 34}/\log X$.   This completes the proof of Theorem \ref{theoremAS}.  
 
 \section{Optimality of Theorem~\ref{theoremAS}} 
  
  \label{optimality} 
  
\noindent In this section, we show the existence of a subset $\Ac$ of $[1,X]$ with $|\Ac| \asymp \sqrt{X}/\log \log X$, and 
such that a positive proportion of the primes below $X$ may be written as $a+ m^2$ with $a\in \Ac$ and $m \in {\Bbb Z}$.  
Since this is only an example to show the optimality of Theorem \ref{theoremAS}, we shall be content with sketching the 
proof.  
  
Set $Z= \exp( (\log X)^{\frac 1{10}})$.  Note that $\log \log Z \asymp \log \log X$.  
Let ${\mathcal D}$ be a set of about $\sqrt{Z}/\log \log X$ odd square-free numbers $d$ with $Z\le d \le 2Z$ and such that   $L(1,\chi_{-4d}) \asymp 1/\log \log X$.   Then our set ${\mathcal A}$ will consist of all numbers of the form $d k^2$ with $d \in {\mathcal D}$ and $k\le \sqrt{X/2Z}$.   By construction the set $\Ac$ has 
$\asymp \sqrt{X}/\log \log X$ elements.

 Arguing using a classical zero-free region for class group $L$-functions, we may see that 
 for any  $d\in {\mathcal D}$ the number of primes up to $X/2$ of the form $dk^2+b^2$ with $b,k \in \N$ is 
  $$ 
  \gg \frac{\pi(X)}{h(-4d)} \asymp \frac{X}{\sqrt{Z}\log X} \log \log X,
  $$
  upon using the class number formula.   
  Thus if $r_{\mathcal A}(p)$ denotes the number of ways of writing $p$ as $a+b^2$ with $a \in \Ac$ and $b \in \N$, it follows that 
  $$ 
  \sum_{p\le X/2} r_{\mathcal A}(p) \gg \frac{X}{\log X}. 
  $$ 
  
By similar methods, we may show that for $d_1 \neq d_2 \in {\mathcal D}$, 
the number of primes up to $X/2$ that may be represented as $d_1 k^2 + b^2$ and also as $d_2 r^2 + s^2$ is at most 
  $$ 
  \ll \frac{\pi( X)}{h(-4d_1) h(-4d_2)} \asymp \frac{X}{Z \log X} (\log \log X)^2. 
  $$ 
  It follows that 
  $$ 
  \sum_{p\le X/2} r_{\mathcal A}(p)^2 \ll \frac{X}{\log X}. 
  $$ 
  By Cauchy-Schwarz it follows that the number of $p\le X/2$ with $r_{\mathcal A}(p) >0$ is  $\gg X/\log X$, as 
  claimed.  
  
In Theorem \ref{theoremAS} we were interested in lower bounds for the size of a set $\mathcal{A}  \subset \{1,\dots, X\}$ such that $\mathcal{A} + \mathcal{B} \supset \mathcal{C}$, where $\mathcal{B}, \mathcal{C} \subset \{1,\dots, X\}$ are given sets (in this case, $\mathcal{B}$ is the set of squares, and $\mathcal{C}$ is a set consisting of a positive proportion of the primes). One might say that $\mathcal{A}$ is an \emph{additive complement of $\mathcal{B}$ relative to $\mathcal{C}$}. In the case $\mathcal{C} = \{1,\dots, X\}$ one recovers the usual notion of additive complement. To our knowledge the relative case has not been studied in any generality. A large number of questions suggest themselves.

\par \bigskip

\end{document}